\documentclass[11pt]{article}
\usepackage{amssymb}
\usepackage{}
\usepackage{pifont}
\usepackage{amscd}
\usepackage{amsfonts}
\usepackage{amsmath, amsthm}
\usepackage{amssymb, amsxtra}

\usepackage{graphicx}

\usepackage{multicol}
\usepackage{epstopdf}
\usepackage{epsf}
\usepackage{float}
\usepackage{extarrows}
\usepackage[a4paper, text={.8\paperwidth,.83\paperheight}, ratio=1:1]{geometry}
\usepackage[format=hang,font=small,textfont=it]{caption}
\usepackage{authblk}

\usepackage[colorlinks=true, linkcolor=blue, citecolor=blue]{hyperref}
\usepackage[english]{babel}
\usepackage{latexsym}

\usepackage{tikz}
\usetikzlibrary{calc}
\usepackage[tight]{subfigure}

\usepackage{color}
\definecolor{greenbean}{RGB}{199,237,204}

\newtheorem{theorem}{theorem}[]

\newtheorem{lemma}{lemma}[]

\newtheorem{example}{example}[]

\newtheorem{definition}{definition}

\newtheorem{remark}{remark}[]

\begin{document}
\title{The Mordell-Weil Groups of Cubic  Pencils
\footnotetext{\hspace{-1.8em} This work is supported by the Natural Science Foundation of Jiangsu Province (No. BK20181427, BK 20211305), and Postgraduate
Research $\&$ Practice Innovation Program of Jiangsu Province~(No. KYCX22$\_$3180)\\
Corresponding author: Mo Jiali E-mail: mojiali0722@126.com\\
.}}

\author[1]{Jia-Li Mo}

\affil[1]{\small{(Department of Mathematics, Soochow University, Suzhou 215006, China)}}

\date{}
\maketitle

\abstract{In this paper we study the influences of the base points of cubic pencils on the Mordell-Weil groups. Specifically, we investigate and classify the cubic pencils with  8, 7 and 6 base points in general position, and give some applications.}

Key words~~~~Base point~~~~Cubic pencil~~~~Elliptic curve~~~~Mordell-Weil group

\section{Introduction}\label{outline}

The elliptic curve is one of the most classical, fundamental and fascinating objects in mathematics. The cubic pencil is an useful method for studying elliptic curves.

Let $H_1,~H_2\in k[X,Y,Z]$~($k$ is an algebraically closed field of characteristic zero)~be two homogenous cubic polynomials without common factor. Consider the cubic pencil
 $$~~~~~~~~~~~~~~~~~~~~~~~~~~~~~~~~~~~~~~~~~~~~~~~~~~~S_{s,t}: sH_1+tH_2=0,~~[s,t]\in \mathbb{P}^1.~~~~~~~~~~~~~~~~~~~~~~~~~~~~~~~~~~~~~~~~~~~~~~~~~~~~~~~~~~~~~~~~~~~~~~~~(1.1)$$
If it contains at least one singular cubic curve, then it defines a genus one fibration over $\mathbb{P}^1$ with homogeneous coordinates $[s,t]$ (possibly after resolving singularities of $S_{s,t}$ as a projective surface in $\mathbb{P}^2 \times \mathbb{P}^1$). Since $k$ is algebraically closed, $S_{s,t}$ gives an elliptic surface with sections given by the base points of the cubic pencil. It is well known that each cubic pencil in $\mathbb{P}^2$ is a rational  elliptic surface by blowing-up of nine (including infinitely near) base points. Moreover, we have the following theorem (~\cite{CD89}, Theorem 5.6.1): {\bf Over an algebraically closed field, every rational elliptic surface
(with section) admits a model as a cubic pencil}.

In  this  paper, we assume that  $S$ is a smooth projective surface having a relatively minimal elliptic fibration $f : S\rightarrow C$ with the zero section $O$ over a curve $C$. Let $E$ be the generic fiber of $f$ which is an elliptic curve over the function field $K= k(C)$. Assume that $f$ has at least one singular fiber. Then the group $E(K)$ of $K$-rational points is finitely generated (Mordell-Weil theorem, see \cite{S2}). It can be identified with the group of sections of $f$. In \cite{Man85},
Manin and Shafarevich computed the Mordell-Weil group of the fibration corresponding to a general cubic pencil in $\mathbb{P}^2$. They proved that the Mordell-Weil group of above fibration is $E_8$.

If we denote by $S$ the blow-up of $\mathbb{P}^2$ at the nine base points, the pencil defines an
elliptic fibration $f: S \rightarrow \mathbb{P}^1$ such that the nine exceptional curves arising from the
blow-up give nine sections $P_0, \cdots, P_8$ of $f$. Choose $P_0$ as the zero section $O$, and
let $E$ be the generic fibre of $f$. Then $P_1, \cdots, P_8$  generate a subgroup of index $3$ in
$E(K)$, and there is an unique point $Q\in E(K)$ such that $\sum_{i=1}^{8}P_i = 3Q$. Together
with $P_1, \cdots, P_7$, the point $Q$ generates the full Mordell-Weil group $E(K)$.

The above result is quite striking since it concerns a general rational elliptic
surface (with section): all cubic pencils with Mordell-Weil rank less than $8$ lie on a
hypersurface inside the moduli space of cubic pencils (the discriminant divisor). It
might therefore come as a little surprise that the general case seems to be exceptional
within rational elliptic surfaces. Indeed, papers \cite{Fus06, Pas12, Salg09} show that
for any cubic pencil of Mordell-Weil rank 4 to 7, the base points generate the full
Mordell-Weil group. At the other end of the scale, \cite{Bea82, Nar87} solved the case
of Mordell-Weil rank zero where there is also a model as a cubic pencil whose base
points generate the finite Mordell-Weil group.

Moreover, in \cite{OS}, Ogusio and Shioda classified the Mordell-Weil groups of a rational elliptic surface.
In  \cite{GKL}, \cite{KK} and  \cite{SSS}, the authors also computed the Mordell-Weil groups of the fibrations corresponding to the curve pencils of  high degree in $\mathbb{P}^2$.

It is well known that 9 points  in  general  position  determine a unique cubic curve in $\mathbb{P}^2$. In this paper, we will give  more accurate results on the Mordell-Weil groups of fibrations corresponding to
the cubic pencils with 8, 7 and 6 base points in general position ($n$ points in general position means that there is no 3 points on a line and no 6 points on a conic).
For the Mordell-Weil groups of such pencils, we have the following theorem.

\begin{theorem}\label{thm1}
Given  $n~~(=8, 7, 6)$ points in general position in $\mathbb{P}^2$. Then in the following cases:

  $(1)$ 8 points are simple base points of a cubic pencil (1.1) (in fact there are 9 simple base points),

  $(2)$ 7 points are only simple base points of a cubic pencil (1.1), and every element of the above cubic pencil is irreducible, and

  $(3)$ 6 points are only simple base points of a cubic pencil (1.1), and every element of the above cubic pencil is irreducible,

\noindent the Mordell-Weil groups of fibrations related to the corresponding cubic pencils are $ E_8,  E^{\vee}_7, E^{\vee}_6$ respectively.
\end{theorem}

\begin{remark}
In our paper, the condition of general position is essential.\\
\begin{example}
In \cite{Bea82}, Beauville gave a cubic pencil (elliptic surface) $$S_{s,t}: ~s(X+Y)(Y+Z)(Z+X)+tXYZ=0,$$ where $H_1=(X+Y)(Y+Z)(Z+X), ~H_2=XYZ$. In this case, there are $8$ base points not in general position~(it has 3 base points on a line ), and rank$(E(K))=0$.
\end{example}
\end{remark}

Moreover, we have the converse form of Theorem \ref{thm1}:
\begin{theorem}\label{thm2}
Let $f : S \rightarrow C=\mathbb{P}^1$ be an elliptic fibration with the zero section $O$. If the Mordell-Weil groups $E(K)\cong E_8$, $E^{\vee}_7$ or $E^{\vee}_6$, then $f$ can be obtained by blowing-up cubic pencils with exactly  9, 7, 6 simple base points respectively.
\end{theorem}

\begin{remark}
In fact, Theorem \ref{thm2} is implied in \cite{S91} in an unobvious form. In this paper,  we provide a new proof for it.
\end{remark}


The cubic pencils of the above Theorem \ref{thm1} can be classified as follows:
\begin{theorem}\label{thm3}
Each of the fibrations with Mordell-Weil groups $E_8,  E^{\vee}_7, E^{\vee}_6$ corresponding respectively to the cubic pencils $(1)-(3)$ in Theorem \ref{thm1} is isomorphic to one of the following two types of fibrations with respect to each Mordell-Weil group:
\begin{equation}
E_8:~~y^2=x^3+x(\sum\limits_{i=0}^3 p_i t^i)+\sum\limits_{i=0}^3 q_i t^i+t^5,~~ y^2=x^3+t^2x^2+x(\sum\limits_{i=0}^2 p_i t^i)+\sum\limits_{i=0}^4 q_i t^i+t^5;\\
\end{equation}
\begin{equation}
E^{\vee}_7:~~y^2=x^3+x(p_0+p_1t+t^3)+\sum\limits_{i=0}^4 q_i t^i,~~ y^2+txy=x^3+x(\sum\limits_{i=0}^2 p_i t^i)+\sum\limits_{i=0}^3 q_i t^i-t^4;\\
\end{equation}
\begin{equation}
E^{\vee}_6:~~y^2+t^2y=x^3+x(\sum\limits_{i=0}^2 p_i t^i)+(\sum\limits_{i=0}^2 q_i t^i),~~ y^2+txy=x^3+x(\sum\limits_{i=0}^2 p_i t^i)+(\sum\limits_{i=0}^3 q_i t^i).\\
\end{equation}
\end{theorem}

A Del Pezzo surface $X$ is either $\mathbb{P}^1  \times \mathbb{P}^1$ or the blow-up of $\mathbb{P}^2$ in $m~(m=1, \cdots ,8)$ points in general position.
The degree $d$ of  $X$  is defined as $d= 9-m$ (see \cite{BHPV}).
As an application, we give a new proof to the numbers of $(-1)$-curves in Del Pezzo surfaces~(see \cite{Man86}).

\begin{theorem}\label{th4}
There are 240, 56, 27~$(-1)$-curves in the Del Pezzo surface of degree $1,2,3$, respectively.
\end{theorem}


The paper is organized as follows. In Section \ref{method}, we recall some related definitions and notations. Then, in Section \ref{proof}, \ref{proof2} and \ref{proof3}, we prove Theorem \ref{thm1}, \ref{thm2} and \ref{thm3} respectively. Finally, in Section \ref{app}, we give some applications.


\section{Definitions and notations}\label{method}

In this section, we recall some related definitions and terminologies.

\begin{definition}
A base point in a cubic pencil is simple if it is a normal crossing point of two general elements in the cubic pencil.
\end{definition}

Now, we recall some basic facts  about the lattices in this paper.
\begin{definition}\label{metho}

A lattice $L$ of rank $r$ is a root lattice of type $E_r$, if there exists a basis $\{\alpha_1, \alpha_2, \cdots, \alpha_r\}$ of $L$ such that for $1\leq i < j\leq r$, we have zero pairing $<\alpha_i, \alpha_j>=0$ unless $(E_r): ~<\alpha_i, \alpha_j>=-1$ for $i+1=j<r$, or ~$i=3, j=r.$
\end{definition}

\begin{figure}[ht]
\begin{center}
\scalebox{0.8}{\includegraphics{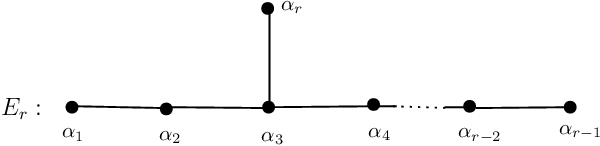}}
\end{center} \caption{Dynkin diagrams of type $E_r$}\label{E}
\end{figure}

In this paper, we are especially concerned with $E_6, E_7$, and $E_8$. For instance, we choose a basis $\{\alpha_1, \alpha_2, \cdots, \alpha_6\}$ of $E_6$ such that  $<\alpha_i, \alpha_i>=2~(i=1, 2, \cdots, 6)$ and ~ $<\alpha_i, \alpha_j>=-1$ for $i+1=j<6$, or $i=3,~j=6$~(see Figure \ref{E}).

\begin{definition}
The dual lattice $L^{\vee}$ of an integral lattice $L$ is defined by $$L^{\vee}=\{x\in L\otimes \mathbb{Q}\mid <x,y>\in \mathbb{Z}, ~\forall~ y\in L\},$$ with pairing naturally extended from $L$ to $L\otimes \mathbb{Q}$.
\end{definition}
As the name suggests, the dual lattice can also be defined as a natural lattice structure on the dual module $L^{\vee}=Hom(L,\mathbb{Z})$. Though this viewpoint is useful in many contexts, we will always regard $L^{\vee}$ as a lattice containing $L$ as a sublattice of finite index.

Now, we recall and fix some standard notation in dealing with Mordell-Weil lattices (cf.~\cite{Shio91}).  The reader can refer to \cite{S2, OS, S1} for more details.
\begin{itemize}
  \item  $k$ : an algebraically closed field of characteristic zero.
  \item  $K=k(C)$ : the function field of $C$ over $k$.
  \item $E(K)$ : the finitely generated Mordell-Weil group, i.e. the group of $K$-rational points of an elliptic curve $E$ over $K$ with the zero section $O$.
  \item  $<P, Q>$ : the height pairing $(P, Q \in E(K))$, as defined by Shioda. it is a symmetric bilinear pairing and a  positive-definite modulo torsion.
  \item  $f : S \rightarrow C=\mathbb{P}^1$ : the associated elliptic surface (the Kodaira-Neron model) of $E/K$ with at least one singular fiber. A $K$-rational point $P\in E(K)$ is identified with a section of $f$.
  \item  $(P)$ : the curve on $S$ is determined by a section $P$, esp. $(O)$  is the zero section viewed as a curve on $S$.
  \item  $(P,O)$ : the intersection number of $(P)$ and $(O)$.
  \item  $\chi(S)$ : the arithmetic genus of $S$ (a positive integer).
  \item  $R~:~= \{v \in C| f^{-1}(v$)~ is~ reducible$\}$.
  \item  $m_v$ : the number of components of fiber $F_v$ of an elliptic surface.
  \item  $\Theta_{v,i}$~$(0\leq i\leq m_{v-1})$ : irreducible components of $f^{-1}(v)(v\in R)$, with $i=0$ corresponding to the identity component, which intersects with the zero section $(O)$.
  \item  $T_v$ : the lattice generated by $\Theta_{v,j}~(j>0)$ with the sign changed. There are root lattices of type $A, D, E$ determined by the type of reducible fiber $f^{-1}(v)$.
 \item $T=\bigoplus_{v\in R}T_v$ : the trivial lattice.
  \item  $E(K)^{0}$ : the narrow Mordell-Weil lattices of $E/K$ defined as $ E(K)^{0}=\{ P\in E(K)|(P)~~{\rm meet}~~\Theta_{v,0},~{\rm for~all}~v \in R \}$. It is a certain subgroup of finite index in $E(K)$.
 \item $\rho(S)$ : the rank of N\"{e}ron-Severi group.
 \item $\equiv$ : numerical equivalence (in the case of elliptic surface, numerical equivalence is equivalent to algebraic equivalence).
 \item $I_p(H_1, H_2)$ : the intersection number of $H_1$ and $H_2$ at $p$.
 \item $l_p(H_1)$ : the multiplicity of curve $H_1$ at $p$. When $p$ is a singular point, $l_p(H_1)\geq 2.$
\end{itemize}

%

\section{Proof of Theorem 1.1}\label{proof}
In this section, we will prove Theorem \ref{thm1}. Let us firstly give several lemmas.

\begin{lemma}\label{lemma1}
Given a rational elliptic surface $f: S\rightarrow \mathbb{P}^1$ with zero section $O$. Let $E$ be a generic fiber; If $C$ is a $(-k)$-curve ~$(k\geq2)$ over $S$, then $C$ is a component of singular fiber and $k=2$.  If $C$ is a $(-1)$-curve, then  $C$ is a section.
\end{lemma}

\begin{proof} Suppose that $C$ is not a component of singular fiber $F$, there would be $(C,F)>0$.
Since we have $\chi(S)=1$, we conclude from $K_S\equiv(2g(\mathbb{P}^1)-2+\chi(S))F$ in \cite{SS}, that $K_S\equiv-F$. It follows that
$C^2=-k,~g(C)=0$ and $(K_S, C)+C^2=2g(C)-2=-2$, then $(K_S, C)\geq0$. This gives $(K_S, C)=-(F, C)\leq0$, this contradicts our assumption. Hence
$(F,C)=0$ and consequently $C$ is a component of singular fiber. From $(K_S, C)=0$ we can get $C^2=-k=-2$.

If $C$ is a $(-1)$-curve, we can get $(K_S, C)=-1$  and $K_S\equiv-E$ by the above two formulas. Thus,  $(K_S, C)=((-E), C)=-1$, $C$ must be a section.
\end{proof}

\begin{lemma}\label{lemma2}
Let $H_1$ and $H_2$ be two cubic curves.

(i) If $I_p(H_1, H_2)=n~ (n=2,3)$, then there exists at least a smooth cubic in the cubic pencil $sH_1+tH_2$.

(ii)  If $H_1$ is smooth, then there exists a smooth cubic $H_3$ in the above cubic pencil, and $H_3\neq H_1$.

\end{lemma}

\begin{proof}
(i) If $l_p(H_1)=m, m\geq2, l_p(H_2)=n, n\geq2$, then $I_p(H_1, H_2)\geq4$. Hence, there exists at least a cubic that is smooth.

(ii) Suppose $H_1$ is smooth. If  $H_2$ is also smooth, then (ii) is true. If  $H_2$ is not a smooth cubic, then we take $H_3=H_1+H_2$ as a new cubic, which is smooth.
\end{proof}

{\bf In the following, we will prove  Theorem \ref{thm1}.}

Proof of Theorem 1.1  (1) 8 points are in general position. In this case, we consider firstly that this cubic pencil contains a reducible cubic $C$. Then $C$ must be split to a conic and a line. On a conic, there are no more than 5 points in their general positions, and at least 3 points remain on a line~(the base points may be the intersection of the line and the conic), which is a contradiction. So every cubic in this cubic pencil is irreducible.

Secondly, since every base point in the cubic pencil is simple, we get 9 simple base points. By blowing-up such 9 points, we can get $f: S\rightarrow \mathbb{P}^1$ with zero section $O$ and 9 $(-1)$-curves, which are all sections. Thus, every fiber is irreducible.

Thirdly, since $S$ is a rational surface, $\rho(S)=10$,
and  all fibers are irreducible, we have rank$(E(K))=8$ by the formula rank$(E(K))=\rho(S)-2-\sum\limits_{v\in R}(m_v-1)$. From the Main Theorem in \cite{OS}, we have $E(K)\cong E_8$.

\begin{figure}[ht]
\begin{center}
\scalebox{0.8}{\includegraphics{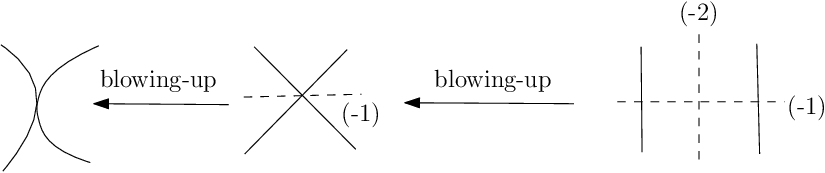}}
\end{center} \caption{The blow-up of base points}\label{E7blowingup}
\end{figure}

 (2) 7 points are in general position will be proved.
In this case, every cubic in this cubic pencil is irreducible. As we have 7 base points in the cubic pencil, which are simple base points, ~according to Lemma \ref{lemma2}, which 2 smooth cubics can be constructed which are tangent at $p$. Then by blowing-up these 7 base points and the multiple base point, we obtain $f: S\rightarrow \mathbb{P}^1$ with zero section $O$, 8 $(-1)$-curves and one $(-2)$-curve ~(see Figure \ref{E7blowingup}). It is clear that $(-1)$-curves must be sections. According to Lemma \ref{lemma1}, the $(-2)$-curve must be a component of  a singular fiber. Hence we just have one singular fiber with 2 components. Then, by the formulas rank$(E(K))=\rho(S)-2-\sum\limits_{v\in R}(m_v-1)$, and $\sum\limits_{v\in R}(m_v-1)=1$, we have rank$(E(K))=7$. By the Main Theorem in \cite{OS}, we have $E(K)\cong E^{\vee}_7$.

(3) 8 points are in general position. The proof is similar to Case (2).

The proof is complete.

\section{Proof of Theorem 1.2}\label{proof2}

\begin{figure}[ht]
\begin{center}
\scalebox{0.8}{\includegraphics{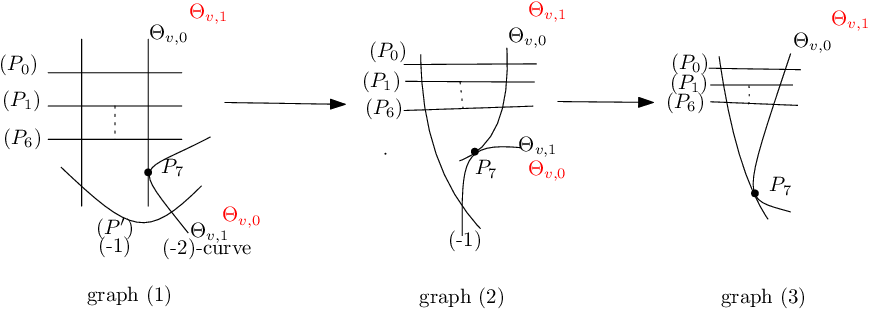}}
\end{center} \caption{The process of contracting $(-1)$-curves}\label{E7}
\end{figure}

{\bf Proof (1) $E(K)\cong E_8$.}
Firstly, we construct 8 $(-1)$-curves in $E(K)$, which are disjoint to each other.
Since $E(K)^{0}\cong E_8$, we only need to consider the above $(-1)$-curves in $E(K)^{0}$. 
Let $(P_0)=O,~ (P_i)=\sum\limits_{n=1}^i \alpha_n ~(i=1,2,\cdots,7)$~(Cf. Definition~\ref{metho} and Figure \ref{E7}). Then  $<P_i,P_i>=2$~and $<P_i,P_j>=1~ (i\neq j, 1\leq i,j\leq7)$. By the formula in \cite{SS}, $<P_i,P_j>=\chi(S)+(P_i, O)+(P_j, O)-(P_i,P_j), ~<P_i,P_i>=2\chi(S)+2(P_i, O)$. We have $(P_i,P_j)=0$ and $(P_i,O)=0$.

Then, we contract the 8 $(-1)$-curves $(P_i) ~(i=0,1,2,\cdots,7)$ to get a cubic pencil $C_t$ with 8 base points. It is well-known that every smooth elliptic curve can be isomorphically mapped to a smooth cubic in $\mathbb{P}^2$. Now, we choose two general curves in  $C_t$, and  isomorphically map them into $\mathbb{P}^2$. By fixing 4 simple base points in $\mathbb{P}^2$, and letting the above 4 points become the images of 4 base points of $C_t$  in $\mathbb{P}^2$, we get a unique cubic pencil in $\mathbb{P}^2$. By the construction, we see that there must be 8 simple base points in the planar cubic pencil. Then by the Bezout theorem, the cubic pencil must have 9 simple base points.

{\bf (2) $E(K)\cong E^{\vee}_7$.}

In this case, we need to construct 7 $(-1)$-curves in $E(k)^0$ which are disjoint to each other. Let $(P_0)=(O), ~(P_i)=\sum\limits_{n=1}^i \alpha_n ~(i=1,2,\cdots,6)$. Similar to  the proof of the case $E(K)\cong E_8$, we have $(P_i,P_j)=0$ and $(P_i,O)=0$.

Then, we construct another $(-1)$-curve $(P')$ such that

(1) $(P')$ is disjoint to $\Theta_{v,0};$

(2) $(P',P_i)=0~(i=0,1,2,\cdots,6).$

Denote the Gram matrix of $E_7$ by $G$, and  $\beta_i,~i=1,2,\cdots,7$  are the basis of $E^{\vee}_7$. Let $G'$ be the
Gram matrix of  $E^{\vee}_7$. Then we have

\small{
\begin{displaymath}
\mathbf{G} =
\left(\begin{array}{ccccccc}
    2 & -1 & 0 & 0 & 0 & 0 & 0 \\
    -1 & 2 & -1 & 0 & 0 & 0 & 0 \\
    0 & -1 & 2 & -1 & 0 & 0 & -1 \\
    0 & 0 & -1 & 2 & -1 & 0 & 0 \\
    0 & 0 & 0 & -1 & 2 & -1 & 0 \\
    0 & 0 & 0 & 0 & -1 & 2 & 0 \\
    0 & 0 & -1 & 0 & 0 & 0 & 2 \\
  \end{array}\right),
  \mathbf{~~~~~G'} =
\left(\begin{array}{ccccccc}
    2 & 3 & 4 & 3 & 2 & 1 & 2 \\
    3 & 6 & 8 & 6 & 4 & 2 & 4 \\
    4 & 8 & 12 & 9 & 6 & 3 & 6 \\
    3 & 6 & 9 & 15/2 & 5 & 5/2 & 9/2 \\
    2 & 4 & 6 & 5 & 4 & 2 & 3 \\
    1 & 2 & 3 & 5/2 & 2 & 3/2 & 3/2 \\
    2 & 4 & 6 & 9/2 & 3 & 3/2 & 7/2 \\
  \end{array}\right).
  \end{displaymath}}

Now take $P'=\beta_6$, since ~~$3/2=<\beta_6,\beta_6>=2\chi(S)+2(\beta_6, O)-\sum\limits_{v\in R} contr_v(\beta_6)$, we have $(\beta_6, O)=0, ~ \sum\limits_{v\in R} contr_v(\beta_6)=1/2$.
Note that $G=GG^{-1}G^{T}$, then $(\alpha_1,\alpha_2,\cdots,\alpha_7)^T=G(\beta_1,\beta_2,\cdots,\beta_7)^T$. So we can use $\beta_1,\beta_2,\cdots,\beta_7$
to represent $\alpha_1,\alpha_2,\cdots,\alpha_7$, hence to represent $P_i, i = 0, 1, 2, \cdots,6$, and then get $<\beta_6, P_i>=1$.  Since $<\beta_6, P_i>=\chi(S)+(\beta_6, O)+(P_i, O)-(\beta_6, P_i)-
\sum\limits_{v\in R} contr_v(\beta_6,P_i)$, and~$\sum\limits_{v\in R} contr_v(\beta_6,P_i)=0$, we get $(\beta_6, P_i)=0 ~(i=1,2,\cdots,6)$.

Since $E(K)\cong E^{\vee}_7,~ T=A_1$, the singular fiber is $\uppercase\expandafter{\romannumeral 3}$ (or $I_2$).
The process of contracting  $(-1)$-curves is as follows: Firstly, we start to blow down $(P')$ and get graph~(2) in Figure \ref{E7}. Then by blowing down $\Theta_{v,1}$ we get graph~(3) in Figure \ref{E7}. Finally we blow down $(P_i)~(i=0, \cdots, 6)$ in turn and get a cubic pencil $C_t$ with 7 simple base points. Then we can also  isomorphically map them into $\mathbb{P}^2$ and fix 4 simple base points in $\mathbb{P}^2$, and get a unique cubic pencil in $\mathbb{P}^2$. The analysis here is similar to the proof for the case $E(K)\cong E_8$.

{\bf (3) $E(K)\cong E^{\vee}_6$}. The proof for this case is similar to the case $E(K)\cong E^{\vee}_7$.

The proof is complete.

\section{Proof of Theorem 1.3}\label{proof3}
\begin{lemma}\label{th}
Let  $f: S\rightarrow \mathbb{P}^1$ be a fibration with the zero section $O$, and $E$ is a generic fiber. If it has no singular fibers, then $S$ can't be a rational surface.
\end{lemma}

\begin{proof}
Assume $f: S\rightarrow \mathbb{P}^1$ is a fibration with no singular fibers. Because there is no Kodaira  fibration over $\mathbb{P}^1$, then $f$ must be a trivial fibration, thus $S=E\times \mathbb{P}^1$ is not a rational surface.
\end{proof}

\begin{lemma}\label{mo}
 If $H_1$ and $H_2$ are two cubic curves in $\mathbb{P}^2$, and $I_p(H_1, H_2)\geq2$. Then there exists a cubic, which  is singular at $p$ in the cubic pencil $sH_1+tH_2$.
\end{lemma}

\begin{proof}
If $H_1$ or $H_2$ is singular at $p$, then the conclusion is obviously true. Now assume $H_1$ and $H_2$ are smooth at $p$. Since $I_p(H_1, H_2)\geq2$,
$H_1$ and $H_2$ have the same tangent at $p$. Then we can choose a number $k$ to make  the coefficient of the lower degree term ($\mathrm{degree} \leqslant1$) of $k H_1-H_2$ to be 0.  Hence, there exists a cubic singular at $p$ in the cubic pencil $sH_1+tH_2$.
\end{proof}

{\bf Proof of Theorem 1.3}

Firstly, when the Mordell-Weil group of fibrations related to the corresponding cubic pencil is $E_8$, we know that  every element in the cubic pencil is irreducible by Theorem \ref{thm1}. From Lemma \ref{th}, we see there is a singular element in the cubic pencil: $sH_1+tH_2$. Because  every cubic is irreducible, then the singular cubic is either a cusp curve or a node curve. In the former case, let $H_0$ be the cusp curve and take it as the form $y^2=x^3$ after some suitable coordinate transforms. Let $H_3$ be another smooth curve, then $\{sH_1+tH_2\}=\{sH_0+tH_3\}$ and $H_0\bigcap H_2=\{p_1, p_2, \cdots, p_9\}$. From the Cayley-Bacharach theorem, $p_9$ can be determined by $p_1, p_2, \cdots, p_8$, and further more,  each $p_i$ can be represented by $v_i$ on $H_0$:
$p_i=(v_i^2,v_i^3), i=1, 2,\cdots,8$.

After blowing-up all the above base points, every $p_i$ corresponds to a section $P_i$ in $E(K)$. Assume $P_9 = O$. Let $u_i=1/v_i$, then $u_i$ corresponds to section  $P_i, ~i=1, 2, \cdots, 8$.  Each integral section $P$  corresponds to a  point $p$ in $H_0$. Because $E(K)\cong E_8$, and rank $\{P_1,...,P_8\}=8$,  $P$ is the unique $\mathbb{Q}$ coefficient combination of $P_i$. Moreover, $p$ is the unique $\mathbb{Q}$ coefficient combination of $\{u_1,\cdots, u_8\}$ by the additive law of $H_0$. Then, every  integral section $P$ corresponds to the unique $\mathbb{Q}$ coefficient combination of $\{u_1,\cdots, u_8\}$, and so the parameters $u_i, i=1, 2, \cdots, 240$ correspond to  240 integral sections. According to the construction theorems in \cite{S3} (Theorem $(E_8)$), the values of $p_i, q_i~(i=0, 1, 2, 3)$ can be obtained from $u_i~(i=1, 2, \cdots, 8)$.

The similar proof can apply to the later case--a node curve (see \cite{S4}).

For the other two cases--the Mordell-Weil groups of fibrations related to the corresponding cubic pencils are $E^{\vee}_7, E^{\vee}_6$, we note that by Lemma \ref{mo}, not all the elements in the cubic pencil are smooth and then the conclusions follow from Theorem $(E_7)$,  Theorem $(E_6)$  in \cite{S4} and \cite{KS}.

\section{Applications}\label{app}

%
%
%
%
%
%
Let us rewrite Theorem \ref{th4} the following form.
\begin{theorem}\label{thm4}

~(1)~ Given 8 points in general position in $\mathbb{P}^2$. If surface $S$ is obtained by blowing-up $p_i, ~i=1,2,\cdots,8$, then there exist exactly 240 $(-1)$-curves on $S$.
 ~(2)~ Given 7 points in general position in $\mathbb{P}^2$. If surface $S$ is obtained by blowing-up $p_i, ~i=1,2,\cdots,7$, then there exist exactly 56 $(-1)$-curves on $S$.
~(3)~ Given 6 points in general position in $\mathbb{P}^2$. If surface $S$ is obtained by blowing-up $p_i,~ i=1,2,\cdots,6$, then there exist exactly 27 $(-1)$-curves on $S$.

\end{theorem}

\begin{proof}
(1) Given 8 points $p_1,p_2,\cdots,p_8$ in general position. It is easy to  construct a cubic pencil $sH_1+tH_2$ with simple base points $p_1,p_2,\cdots,p_8$.
Then, by blowing-up these base points, we can get a surface $\tilde{S}$ and 9 sections $P_i~i=1,2,\cdots,9$. By Theorem \ref{thm1} $E(K)\cong E_8$.
So the number of integral sections is 240.

Assume $(P_9)=(O)$. Then by blowing-down $(P_9)$ on $\tilde{S}$, we obtain the surface $S$. Note that the
$(-1)$-curves in $S$. result from the $(-1)$-curves or  $(-2)$-curves on $\tilde{S}$. If a $(-1)$-curve on $S$ results from a  $(-1)$-curves on $\tilde{S}$, then this  $(-1)$-curve can't intersect with $(O)$ on $\tilde{S}$, and so it is an integral section on  $\tilde{S}$. Hence the number if $(-1)$-curve on $S$ is 240. On the other hand, by Lemma \ref{lemma1}, the $(-2)$-curves on $\tilde{S}$ must be irreducible components of reduced fibers. Since all fibers are irreducible, the $(-2)$-curves never exist on $\tilde{S}$.
\end{proof}

%
%

 \begin{proof}
(2) Firstly, according to the construction in \cite{Shio95a}, there exists a cubic pencil $sH_1+tH_2$ such that $H_1$ and $H_2$ which intersect properly at $p_0, p_1, \cdots, p_6$ and at another 2 multiple point $p_7$.

Then, by blowing-up all the base points of $sH_1+tH_2$, we can get a surface $\tilde{S}$. The only reduced singular fiber is $\uppercase\expandafter{\romannumeral 3}$ (or $I_2$), and $E(K)=E^{\vee}_7$  by Theorem \ref{thm1}. Now assume $(P')=(O)$. Since $\Theta_{v,0}$ corresponds to the identity component, which intersects with the zero section $(O)$, (see Figure \ref{E7} and  the red fiber components there),  by blowing down $(P')$ ans $\Theta_{v,1}$ in turn, we obtain $S$.

Note that the $(-1)$-curves on $S$, result from the $(-k)$-curves on $\tilde{S}$. We obtain the number of $(1)$-curves on $S$ as follows.

Since a $(-1)$-curve must be a section (see Lemma \ref{lemma1}), it just intersects with one of singular fiber components. If a $(-1)$-curve on $S$ intersects with $(O)$ or $\Theta_{v,0}$,  then intersection number  will go  up after contraction. If the $(-1)$-curves on  $S$ come from  the $(-1)$-curves on $\tilde{S}$, then  the $(-1)$-curves on $S$ must intersect with $\Theta_{v,1}$. Since the norm of sections $<P,P>=2\chi(S)+2(O,P)-\sum\limits_{v\in R} contr_v(P)=3/2$, the number of such sections is 56 (see \cite{S2}).

Now consider the $(-1)$-curves on $S$ result from the $(-k)$-curves $(k\geq2)$ on $\tilde{S}$. Since a  $(-k)$-curve must be one component of a reducible singular fiber ( see Lemma \ref{lemma1}), and $\Theta_{v,0}$ is contracted on $S$, after two times of contracting, the self intersection number of  $\Theta_{v,1}$ is a positive integral number. Hence, in this case, there is no $(-1)$-curves resulting from $(-k)$-curves.

In conclusion, the number of $(-1)$-curves on $S$ is 56.
\end{proof}

The number of $(-1)$-curves on $S$ for case (3) can be obtained similarly to case (2). Readers can refer to \cite{Fus06} and \cite{Shio95a}.

\begin{remark}
We can also can get some special fibrations for the cases (1) and (2) in Theorem \ref{thm3}. When $p_i=0, q_j=0$ for all $i$ and $j$, the first equation in (1.2) becomes $y^2=x^3+t^5$, which is a fibration over $\mathbb{P}^1$ with two singular fibers. When all $p_i=0$, and $q_0=q_1=q_2=q_3=0,~q_4=1$, the first equation (1.3)  becomes $y^2=x^3+xt^3+t^4$, which is a fibration over $\mathbb{P}^1$ with three singular fibers. For more details see \cite{SH}.
\end{remark}

\end{document}